\documentclass[12pt]{amsart}
\usepackage{amsmath}
\usepackage{amssymb}
\newcommand{\vp}{\varphi}
\newcommand{\begeq}{\begin{equation}}

\newcommand{\zzendeq}{\end{equation}}
\newcommand{\begmat}{\left(\begin{array}{cc}}
\newcommand{\zzendmat}{\end{array}\right)}
\newcommand{\R}{\mathbf{R}}

\newcommand{\Q}{\mathbf{Q}}

\newcommand{\N}{\mathbf{N}}

\newcommand{\Sum}{\sum\limits}

\newcommand{\la}{\langle}
\newcommand{\ran}{\rangle}

\newcommand{\tr}{\text{tr}}
\newcommand{\dra}{\rangle\rangle}
\newcommand{\dla}{\langle\langle}

\newtheorem{theorem}{Theorem}[section]
\newtheorem{lemma}[theorem]{Lemma}

\newtheorem{proposition}[theorem]{Proposition}

\numberwithin{equation}{section}
\newcommand{\abs}[1]{\lvert#1\rvert}

\begin{document}
\title[Equidistribution of measures on the sphere]{Arithmetic and equidistribution
  of measures on the sphere} 
\author[S. B\"ocherer \and  P.  Sarnak \and
  R. Schulze-Pillot]{Siegfried B\"ocherer  \and  Peter Sarnak \and
  Rainer Schulze-Pillot}  
\thanks{Part of this work was done during a stay of all three authors
  at the IAS, Princeton, supported by von Neumann Fund, Weyl fund and Bankers
  Trust Fund. We thank the
  IAS for its hospitality.}  
\thanks {S. B\"ocherer and R. Schulze-Pillot thank D. Prasad and the Harish Chandra
  Research Institute, Allahabad, India for their
  hospitality. Schulze-Pillot's visit to HCRI was also supported by DFG}  

\maketitle
\begin{abstract}
Motivated by problems of mathematical physics (quantum
chaos) questions of equidistribution of eigenfunctions of the Laplace 
operator on a Riemannian manifold have been studied by several
authors. We consider here, in analogy with arithmetic hyperbolic
surfaces, orthonormal bases of eigenfunctions of the Laplace operator
on the two dimensional unit sphere which are also eigenfunctions of
an algebra of Hecke operators which act on these spherical harmonics.
We formulate an analogue of the equidistribution of mass conjecture
for these eigenfunctions as well as of the conjecture that their moments
tend to moments of the Gaussian as the eigenvalue increases. For such
orthonormal bases we show that these conjectures are related to the
analytic properties of degree eight arithmetic L-functions associated
to triples of eigenfunctions. Moreover we establish the conjecture for 
the third moments and give a conditional (on standard analytic conjectures
about these arithmetic L-functions) proof of the equdistribution of
mass conjecture.

\end{abstract}
\section{Introduction}
\label{intro}
Let $X$ be a Riemannian manifold of finite volume. Starting out from
problems of theoretical physics (quantum chaos) several authors have
recently studied questions of equidistribution of  eigenfunctions of
the Laplace operator. 

In particular, 
precise versions of conjectures on equidistribution properties 
have been put forward by
Rudnick and Sarnak \cite{Ru-Sa} for arithmetic hyperbolic manifolds
$X_\Gamma=\Gamma\backslash H,$ where $H$ is the upper half plane of the complex
numbers and $\Gamma$ an arithmetic subgroup of $SL_2(\mathbf R)$ and
for eigenfunctions of the Laplace operator, that are eigenfunctions of
the (arithmetically defined) Hecke operators as well. 
The phenomenon of (conjectural)
equidistribution of eigenfunctions in this arithmetical situation is
one of the central problems in what has
become known as arithmetic quantum chaos.

We investigate here the analogous question for the
situation of the $2$-dimen\-sio\-nal unit sphere. Although the
dynamics of geodesics for this manifold is certainly not chaotic it turns out that
it nevertheless makes sense to look for an equidistribution
property of eigenfunctions. At first sight, the well known fact that
the usual spherical eigenfunctions $Y_{l,m}$ (see \cite[Chapter
III]{vilenkin}) concentrate for $l=m
\rightarrow \infty$ around the equator \cite{colindeverdiere}seems to contradict the expectation of
equidistribution, but since the eigenvalues occur on the sphere with
multiplicities bigger than one, it makes sense to look into the
question what happens  if one varies the basis of eigenfunctions.

In this direction, it has been proved by Zelditch \cite{zelditch} that
for a random orthonormal basis of eigenfunctions the equidistribution
of mass conjecture is true.

We consider here, in analogy to the arithmetic hyperbolic surfaces,
an orthonormal basis of eigenfunctions of the Laplace operator that
are also eigenfunctions of an algebra of Hecke operators that acts on
the space of spherical functions.
The papers \cite{jakobson-zelditch} and \cite{vanderkam}
also consider questions of the behaviour of
 eigenfunctions for such bases of spherical harmonics. 

Concretely, a definite quaternion algebra such as the Hamilton
Quaternions  $\mathbb H$ 
over $\Bbb Q$, gives rise to Hecke operators on $L^2(S^2)$
(see \cite{eichler72}, \cite{LPS}). For 
 $$\alpha = x_0+x_1i+x_2j+x_3k \in \Bbb H(\Bbb R),$$
let 
 $$S(\alpha) = \frac{1}{\sqrt{N(\alpha)}} \left[ \begin{array}{rl}
   x_0+x_1i & x_2+x_3i \\
   -x_2+x_3i & x_0-x_1i \end{array}\right] \in SU(2).$$
Here
 $$N(\alpha) = \alpha \overline{\alpha} = x_0^2+x_1^2+x_2^2+x_3^2.$$
For $n \geq 1$ an odd integer define the Hecke operator $T_n$ on
$L^2(S^2)$ by
 $$(T_n\phi)(P) = \sum_{N(\alpha)=n\atop 
  \alpha \in \Bbb H(\Bbb Z)} \phi(S(\alpha)P),$$
where $P \in S^2$ and $SU(2)$ acts on $S^2$ by isometries after one realizes
$S^2$ as $\Bbb C \cup \{\infty\}$ via stereographic projection and
$SU(2)$ acts by linear fractional transformations. The $T_n$'s are 
selfadjoint, they commute with each other as well as with the Laplacian 
$\Delta$ on the round sphere. Thus the $T_n$'s can be simultaneously
diagonalized in each of the $2\nu+1$ dimensional spaces $H_{\nu},$
consisting of spherical harmonics on $S^2$ of degree $\nu$ (that 
is the restriction of harmonic polynomials in $\Bbb R^3$, homogeneous
of degree $\nu$). 
This algebra of Hecke operators arises naturally if one views the spherical
harmonics as components at infinity of automorphic forms on the
multiplicative group of the adelization of the rational Hamilton
quaternions. 
We denote by $\psi_{\nu}$ such a Hecke eigenform with
$\nu$ indicating its degree (so that its Laplace eigenvalue is 
$\nu(\nu+1)/2)$. 

The analogue of the equidistribution of mass conjecture
\cite{Ru-Sa} for the $\psi_{\nu}$'s is the following:
 \vspace{0.3cm}\\
{\bf Conjecture 1.} {\it Normalize $\psi_{\nu}$ on $S^2$ to have
$L^2$-norm equal to $1$, so that}
 $$\mu_{\psi_{\nu}} := |\psi_{\nu}(P)|^2 dv(P)$$
{\it is a probability measure. Then}
 $$\lim_{\nu\to \infty} \mu_{\psi_{\nu}} = \frac{dv}{2\pi},$$
{\it in the sense of integration against continuous functions on $S^2$.}
 \vspace{0.3cm}\\
The analogue of the Gaussian equidistribution conjecture of Berry and
others \cite{hejhal} in this context is as follows:
 \vspace{0.3cm}\\
{\bf Conjecture 2.} {\it Fix $q \geq 0$ an integer then}
 $$\lim_{\nu \to \infty} \int_{S^2} \psi_{\nu}^q dv 
  \longrightarrow \frac{c_q}{(2 \pi)^{q/2}},$$
{\it where $c_q$ is the $q$-th moment of the Gaussian distribution.}
 \vspace{0.3cm}\\

By the work of Eichler \cite{eichler72} and of Jacquet/Langlands it is
known that there is a correspondence between spherical harmonic
polynomials and modular forms via the theory of theta series with
spherical harmonics. This correspondence is Hecke-equivariant, and
thus methods and results from the theory of modular forms, in
particular from the theory of $L$-functions associated to Hecke
eigenforms (or irreducible automorphic representations),  can be
used in the study of the spherical harmonics.

The crucial point for our study of the integrals appearing in the
equidistribution conjecture above is a formula proved in
\cite{bs-triple} that connects the integral of a product of 3
eigenfunctions over the sphere with the central critical value of the
automorphic $L$-functions associated to a triple of modular Hecke
eigenforms; this allows one to connect the equidistribution conjecture
with conjectural properties of such automorphic $L$-functions. We note
in passing that such integrals of products of $3$ eigenfunctions of
the Laplace operator on the sphere have been considered in various
places in the physics literature, see \cite{sebilleau}.

The purpose of this note is to show that combining the main formula in
\cite{bs-triple} with the recent subconvex estimates for special
values of $L$-functions of holomorphic modular forms
\cite{peng} allows one to prove 
Conjecture 2 
for $q=3$ (the cases $q=1$ and $q=2$ are obvious). We also show that 
Conjecture 1 would follow from subconvex estimates for the
degree 8 $L$-functions mentioned above. Such subconvex estimates are
an immediate consequence of the Riemann Hypothesis for these 
$L$-functions. At the present time such subconvex estimates are known only
for special forms, see \cite{SAR} and \cite{IS}. 

In his recent thesis \cite{WAT}, Watson has derived general explicit 
identities relating integrals of products of 3 Maass (or holomorphic)
Hecke eigenforms on arithmetic surfaces, to special values of degree 8
$L$-functions. As a consequence he obtains similar results for
``chaotic'' eigenstates.

\medskip
As an appendix to this paper we give a list of corrections to the
article \cite{bs-triple}, on whose results the estimates in the
present paper depend. A revised version of that article is available
at\\ 
{\tt www.math.uni-sb.de/{\~{}}ag-schulze/Preprints}.

We would like to thank T. Ibukiyama for the permission to use his
unpublished results in \cite{ibu}.

\section{Equidistribution}
\label{equi}
Our first goal is to describe explicitly the connection
between the central critical value of the triple product
$L$-function associated to a triple of cusp forms on one
side and integrals of harmonic polynomials over the unit
sphere on the other side. We fix first some notations.

We consider a definite quaternion algebra $D$ 
of discriminant $N$
(where $N$ is the product of the primes ramified in $D$)
over $\Q$ and 
a maximal order $R$ in $D$, we assume that the class number
(i.e., the number of classes of left $R$-ideals) of $D$ is
1;
this restricts $D$ 
to be one of the algebras of discriminant equal to $2, 3, 5, 7, 13.$

On $D$ we have the involution $x\mapsto\overline x$, the
(reduced) trace
tr$(x)=x+\overline x$ and the (reduced) norm $n(x)=x\overline x$.

\medskip
For $\nu \in \N$ let $U_{\nu}^{(0)} $ be the space of homogeneous
harmonic polynomials of degree $\nu$ on $\R^3$ and view $P \in 
U_{\nu}^{(0)}$ as a polynomial on $D_\infty^{(0)}= \{ x \in
D_\infty =D\otimes \R
\vert \tr (x)=0\}$ by putting
$P(\sum_{i=1}^{3}x_ie_i)=P(x_1,x_2,x_3)$ 
for an orthonormal basis $\{ e_i\}$ of $D_\infty^{(0)}$ with
respect to the norm form $n$. Integrating the polynomial in
this identification over the set of $x \in D_\infty^{(0)}$ of
norm $1$ is the same as integrating the original polynomial
in 3 real variables over the unit sphere $S^2,$ we will
freely use this identification below.

In the same way we fix an orthonormal
basis of $D_\infty$ extending the one from above and use it to
identify (harmonic) polynomials in 4 variables with (harmonic)
polynomial functions on $D_\infty$. 

The representation of $D_\infty^\times$ on $U_{\nu}^{(0)}$ by
conjugation of the argument is denoted by $\tau_\nu$.
By $\dla\quad,\quad\dra$ we denote the 
invariant scalar product in the representation space
$U_\nu^{(0)}$ (where the choice of normalization will be
discussed later). 
The
$D^\times \times D^\times$-space
$U_{\nu}^{(0)}\otimes 
U_{\nu}^{(0)}$
is isomorphic to the $D^\times \times D^\times$-space
$U_{2\nu}$ of harmonic  
polynomials on $D_{\infty}$ of degree $2\nu,$ where
$(d_1,d_2)\in D^\times \times D^\times$ acts by sending
$P(x)$ to $P(d_1^{-1}xd_2).$  
An explicit isomorphism is given by mapping $P_1\otimes P_2$
to the polynomial $P_1 \otimes P_2(d):=\dla P_2(x),P_1(dx\bar{d})\dra.$
We will henceforth identify $U_{\nu}^{(0)}\otimes 
U_{\nu}^{(0)}$ with $U_{2\nu}$ using this isomorphism.

There is a Hecke action on $U_{\nu}^{(0)}$ which has been
described by Eichler \cite{eichler72} in terms of Brandt
matrices with polynomial entries, it is given by
$$\tilde{T}(p)P=\sum_{ y\in R,n(y)=p}\tau_\nu(y)(P),$$
see also \cite{LPS}. In particular the space $U_{\nu}^{(0)}$
has a basis consisting of eigenforms of all the
$\tilde{T}(p)$ for the $p\nmid N.$

To $P_1 \in 
U_{\nu}^{(0)}$ we associate the theta series of $R$ with
harmonic polynomial $P_1\otimes P_1$ given as usual as
$$f_{P_1}(z):=\frac{1}{\abs{R^\times}}\sum_{x\in
  R}(P_1\otimes P_1)(x)\exp(\pi i 
n(x)z).$$ For this to be nonzero we have to restrict to
polynomials $P_1$ that are invariant under the action of the
group $R^\times$ of $R,$ we will always do so in the
sequel. The function $f_{P_1}$ is then a cusp form for
$\Gamma_0(N)$ of weight 
$2+2\nu$ if $\nu>0$ and it is an 
eigenform for the Hecke operators $T(p)$ for $p\nmid N$ if $P_1$ is an
  eigenfunction of the $\tilde{T}(p)$ for the $p\nmid N.$
In
fact it is a normalized newform if $\dla P_1,P_1 \dra =1, $
and it is a result of \cite{eichler72} that one gets all
normalized newforms of level $N$, weight 
$2+2\nu$ and trivial character in this way (we will actually
not use the latter fact).
 With these notations we can now formulate:
\begin{proposition}\label{firstproposition}
Let $P_1,P_2,P_3\in U_{\nu_1}^{(0)},U_{\nu_2}^{(0)},
U_{\nu_3}^{(0)}$ (with $\nu_1=\nu_2,\nu_3>0$) be harmonic
polynomials that are Hecke eigenforms as above, denote by 
$f_1,f_2,f_3$ the associated cusp forms of weights $k_1=k_2=2+2\nu_1,k_3=2+2\nu_3$ and by
$L(f_1,f_2,f_3;s)$ the triple product $L$-function
associated to $f_1,f_2,f_3,$ (as defined for the good primes e.g. in \cite{Ga},
for the Euler factors at the bad primes we refer to
\cite{bs-triple}).\\
Then one has for all $\epsilon>0:$
\begin{equation}
  \label{eq:mainresult1}
  L(f_1,f_2,f_3;\frac{2k_1+k_3}{2}-1)
\ge C_1(N,\nu_3)\nu_1^{1-\epsilon}
  (\int_{S^2} P_1(x)P_2(x)P_3(x)dx)^2
\end{equation}
with a poitive constant $C_1(N,\nu_3)$ depending only on $N,\nu_3$ and $\epsilon.$

If $\nu_1=\nu_2=\nu_3=:\nu,$ one has for all $\epsilon>0$ (with $k:=k_1=k_2=k_3=2+2\nu$):
\begin{equation}
  \label{eq:mainresult1a} L(f_1,f_2,f_3;2+3\nu)
\ge C_2(N)\nu^{2-\epsilon}
  (\int_{S^2} P_1(x)P_2(x)P_3(x)dx)^2
\end{equation}
with a positive constant $C_2 (N)$ depending only on $N$ and $\epsilon.$

\end{proposition}

{\it Proof.}
According to \cite{bs-triple} the central critical value is:

\begin{multline}\label{theformula}
(-1)^{a'}2^{5+4a+3b-\omega(N)}\pi^{5+9a'+4b}
\frac{(a'+1)^{[b]}}{2^{[a+b]}2^{[a']}(b+1)^{[a']}(1)^{[a']}}\\
\times\la f_1,f_1\ran\la f_2,f_2\ran \la
f_3,f_3 \ran\bigl(T_0\bigl(
\frac{1}{\abs{R^\times}}
P_1\otimes
P_2\otimes
P_3\bigr)\bigr)^2 
\end{multline}
with a certain trilinear form $T_0$ on $U_{nu_1}^{(0)}\otimes
U_{nu_2}^{(0)}\otimes U_{nu_3}^{(0)}$ whose definition is recorded below.

Here we have the following notations:

$$ \alpha^{[\nu]}= {\Gamma(\alpha+\nu)\over \Gamma(\alpha)}=
\left\{
\begin{array}{lcl}
1 && \nu=0 \\
\alpha (\alpha +1)\dots (\alpha + \nu -1) & & \nu >0
\end{array} \right\}, $$
hence
$$\frac{(a'+1)^{[b]}}{2^{[a+b]}2^{[a']}(b+1)^{[a']}(1)^{[a']}}=\frac{b!}{(a'!)^2(a'+1)!(a+b+1)!}.$$

\medskip
Unfortunately, \cite{bs-triple} contains a mistake at this point,
the correct value of the factor arising here is

\begin{equation}
  \label{eq:triplecorrect}
  \frac{(a'+1+b)b!}{(a'!)^2(a'+1)!(3a'+b+1)!}.
\end{equation}

Moreover, there should be an additional factor of $N^{-1}$
in  (\ref{theformula}), the exponent at $\pi$ should be
$5+6a'+2b$ and the factor $(-1)^{a'}$ should be omitted.

\medskip
The forms $f_1,f_2,f_3$ are normalized newforms of weights $k_1\ge k_2\ge k_3$
with $k_2+k_3 >k_1, $ in our case we have $k_1=k_2$. We write
$k_1=k_2=2+a+b$ and $k_3=2+a, a=2a'.$ For our purposes we can restrict
to the case that both $a'$ and $b$ are even.

We   normalize the invariant scalar 
product on the latter space in such a way that 
the Gegenbauer
polynomial 
$G^{(\alpha)}(x,x')$ obtained from 
$$G^{(\alpha)}_{1}
(t) \ = \ 2^{\alpha }
\ \mathop\sum\limits^{[{\alpha\over 2}]}_{j=0} (-1)^j \ 
{{1}\over{j!(\alpha- 2j)!}} \ 
{{(\alpha-j)!}\over 2^{2j}} \ t^{\alpha-2j}$$ 
by
$$ G^{(\alpha)}(x_1,x_2)= 2^\alpha
(n(x_1)n(x_2))^{\alpha/2}G_1^{(\alpha)}({{{\tr} (x_1  
\overline{x_2})}\over {2\sqrt{n(x_1)n(x_2)}}})$$
 
is a reproducing kernel.  
 
 The invariant scalar product on
$U_\alpha^{(0)}$ is then normalized such that we obtain the product on
$U_{2\alpha}$ given above under the identification
$U_{2\alpha}=U_\alpha^{(0)}\otimes U_\alpha^{(0)}.$ 

Given this, the polynomials $P_1,P_2,P_3$ 
are normalized to 
$$\dla P_i,P_i\dra=1.$$

The normalization of the trilinear form is then as follows:\\
We have a harmonic polynomial $P\in U_{a+b}\otimes U_{a+b}\otimes
U_{a}$ in three vector variables (each 
vector being a quaternion) derived from the action of a certain
differential operator on an exponential in Section 1 of \cite{bs-triple}.
This gives an invariant trilinear form $T$ on   $ U_{a+b}\otimes
U_{a+b}\otimes U_{a}$ defined by taking the scalar product with
$P_0:=(\pi i)^{-3a'-b}i^{-3a'}P$  (notice that in \cite{bs-triple} we write
erroneously $\pi^{-3a-2b}P$).
Using the identification
$U_{2\alpha}=U_\alpha^{(0)}\otimes U_\alpha^{(0)}$  from above $T$
decomposes as $T_0\otimes T_0.$

\medskip
For the intended application the form $T_0(Q_1^{(0)},Q_2^{(0)},Q_3^{(0)})$ should be
replaced by the integral 
$$ \int Q_1^{(0)}({\bf x})Q_2^{(0)}({\bf x})Q_2^{(0)}({\bf x}) d{\bf x}$$ over the unit sphere.

\medskip
As a first step we 
compare $T(Q_1,Q_2,Q_3)$ with $$ \int Q_1({x})Q_2({x})Q_2({x}) d{
  x};$$
since both expressions give invariant trilinear forms they have to be
proportional. 

We compute $T$ for special polynomials $Q_i$ on the
space of quaternions:\\
Write $G^{(\alpha)}_{ w}({ x})$ for
$G^{(\alpha)}({ w},{ x}).$

Then we have 
$$T(G_w^{(a+b)},G_w^{(a+b)},G_w^{(a)})=P_0(w,w,w)$$
by the reproducing property of the $G^{(\alpha)}.$

On the other hand, by \cite[p.490]{vilenkin} the integral 
$$\int G_w^{(a+b)}(x)G_w^{(a+b)}(x)G_w^{(a)}(x)dx$$ is equal to
$\pi/2$ and hence 
\begin{equation}
  \label{eq:tcompare}
   T(Q_1,Q_2,Q_3)=2P_0(w,w,w)\int Q_1({\bf x})Q_2({\bf x})Q_2({\bf x}) d{\bf
  x}.
\end{equation}

We have to compute $P_0(w,w,w)$ explicitly. This looks at
first sight rather awkward since our description in
\cite{bs-triple} gives us an explicit formula only for one
coefficient of the polynomial. 

Fortunately there are some results on such
polynomials in forthcoming work of Ibukiyama and Zagier, see  
\cite{ibu}:\\
For $n\in {\bold N}$ we denote by ${\mathcal H}_n(4)$ the space of
harmonic homogeneous polynomials of degree $n$ in $4$ variables.
For nonnegative integers $\mu_1,\mu_2,\mu_3$ we then put
$${\mathcal H}_{\mu_1,\mu_2,\mu_3}(4):=
\left( {\mathcal H}_{\mu_2+\mu_3}(4)\otimes 
{\mathcal H}_{\mu_1+\mu_3}(4)\otimes {\mathcal H}_{\mu_1+\mu_2}(4)\right)^{O(4)}$$
This space is then always one-dimensional and a nonzero element of 
${\mathcal H}_{\mu_1,\mu_2,\mu_3}(4)$ is (explicitly!) given as the
 coefficient
of $X_1^{\mu_1}X_2^{\mu_2}X_3^{\mu_3}$ in the formal power series 
$$G_4({\bf X},T)=G_4(X_1,X_2,X_3;T)=
{1\over
\sqrt{\Delta({\bf X},T)^2-4d(T)X_1X_2X_3}}$$
Here $T$ is (twice of) a Gram matrix
$$T=\left(\begin{array}{ccc}
2m_1 & r_3 & r_2\\
r_3 & 2m_2 & r_1\\
r_2 & r_1 & 2m_3 \end{array}\right)=
2\cdot \left(\begin{array}{ccc}
({\bf x},{\bf x}) & ({\bf x},{\bf y}) & ({\bf x},{\bf z})\\
({\bf y},{\bf x}) & ({\bf y},{\bf y}) & ({\bf y},{\bf z})\\
({\bf z},{\bf x}) & ({\bf z},{\bf y}) & ({\bf z},{\bf z})
\end{array}\right)\quad ({\bf x},{\bf y},{\bf z}\in {\bold{C}}^4)$$
and
$$d(T)=4m_1m_2m_3-m_1r_1^2-m_2r_2^2-m_3r_3^2+r_1r_2r_3=
{1\over 2}\det(T)$$
\begin{eqnarray*}
\Delta(X_1,X_2,X_3;T)& = & \Delta({\bf X},T)\\
{} & = & 
1-r_1X_1-r_2X_2-r_3X_3+r_1m_1X_2X_3+
r_2m_2X_3X_1 +r_3m_3X_1X_2\\
{} & {} & +m_1m_2X_3^2+m_2m_3X_1^2+m_3m_1X_2^2.
\end{eqnarray*}

We are interested in the coefficient of $X_1^{a'}X_2^{a'}X_3^{a'+b}$
which we call $\tilde{P}$ in the sequel.
The coefficient of $m_1^0m_2^0m_3^0r_1^{a'}r_2^{a'}r_3^{a'+b}$ in the
polynomial $\tilde{P}$ can be read off from the expression above 
(putting $m_1=m_2=m_3=0$); it is
$$\sum_{\alpha=0}^{a'}{2a'-2\alpha \choose a'-\alpha} {2a'+b+\alpha
  \choose 3\alpha +b}{3\alpha+b\choose \alpha,\alpha,\alpha+b},$$
where we write 
$${j\choose \alpha,\beta,\gamma}:={j!\over \alpha!\beta!\gamma!}.$$
This is known to be equal to
$$\frac{(2a')!(b+2a')!^2}{a'!^4(b+a')!^2};$$
an identity which can be reduced to a special case of an exercise on
page 44 in \cite{stanley} 
with hints to \cite{andrews}, who traces it back to "Saalschutz summation".

Now we compare $P_0$ with Ibukiyama's polynomial $\tilde{P}$;
it is enough to compare the coefficients in the monomial above.
From Section 1 of \cite{bs-triple} one reads off that
the coefficient of $P_0$ in the same monomial is
$$\frac{2^b}{(b+1)!b!}.$$

Again, there is a mistake in \cite{bs-triple} here.
The correct value is
\begin{equation}
  \label{eq:triplecorrect2}
  \frac{2^b 2^{4a'}\Gamma(a'+2)}{\Gamma(a'+b+2)b!},
\end{equation}

so that we arrive at 
$$P_0=2^{b+4a'}\frac{a'!^4(b+a')!^2(a'+1)!}{(2a')!(b+2a')!^2(a'+b+1)!b!}\tilde{P}.$$

Next we have to evaluate $\tilde{P}$ at the triple $(w,w,w)$, i.\ e.,
at the matrix $T$ with $m_1=m_2=m_3=1, r_1=r_2=r_3=2.$
We get in this case  
$$ G_4({\bf X},T)=(1-\sum_{i=1}^3X_i)^{-2},$$
the coefficient of which at $X_1^{a'}X_2^{a'}X_3^{a'+b}$ is the value
we try to compute.

It is proved easily (Taylor expansion) that  
this is equal to 
$$\frac{(3a'+b+1)!}{(a')!(a')!(a'+b)!},$$
which leads us to 

\begin{eqnarray*}P_0(w,w,w)&=&2^{b+4a'} \frac{(a'+1)!(3a'+b+1)!a'!^4(b+a')!^2}{a'!^2 (a'+b)!(2a')!(b+2a')!^2(a'+b+1)!b!}\\
&=&2^{b+4a'} \frac{(3a'+b+1)!a'!^2(a'+1)!(a'+b)!}{(2a')!(b+2a')!^2b!(a'+b+1)!}.
\end{eqnarray*}

\medskip
For polynomials $Q_1^{(0)},Q_2^{(0)}\in U_{a'+b/2}^{(0)},Q_3^{(0)}\in
U_{a'}^{(0)} $ we have by definition
$$(T_0(Q_1^{(0)},Q_2^{(0)},Q_3^{(0)}))^2=T(Q_1^{(0)}\otimes Q_1^{(0)},
Q_2^{(0)}\otimes Q_2^{(0)}, Q_3^{(0)}\otimes Q_3^{(0)})$$ and hence (as
a consequence of the discussion given above)

\begin{multline*}(T_0(Q_1^{(0)},Q_2^{(0)},Q_3^{(0)}))^2=2^{b+4a'+1}\pi^{-1}
  \frac{(3a'+b+1)! a'!^2(a'+1)!(a'+b)!}{(2a')!(b+2a')!^2b!(a'+b+1)!}\\
\times \int
(Q_1^{(0)}\otimes Q_1^{(0)})(x)( Q_2^{(0)}\otimes
Q_2^{(0)})(x)(Q_3^{(0)}\otimes Q_3^{(0)})(x)dx,
\end{multline*}
where the integration is over the 3-dimensional unit sphere. 

\medskip
Our next task is to relate the integral 
\begin{equation}
  \label{eq:int_in4}
  \int
(Q_1^{(0)}\otimes Q_1^{(0)})(x)( Q_2^{(0)}\otimes
Q_2^{(0)})(x)(Q_3^{(0)}\otimes Q_3^{(0)})(x)dx
\end{equation}
with 
$(\tilde{T}_0(Q_1^{(0)},Q_2^{(0)},Q_3^{(0)})^2,$
where we put 
\begin{equation}
  \label{eq:int_in3}
  \tilde{T}_0(Q_1^{(0)},Q_2^{(0)},Q_3^{(0)}):=\int
Q_1^{(0)}(z)Q_2^{(0)}(z)Q_3^{(0)}(z)dz,
\end{equation}

and where the integration is now over the 2-dimensional unit sphere
(the factor of proportionality arising here depends on the
identification between $U_{2\alpha}$ and $U_\alpha^{(0)}\otimes
U_\alpha^{(0)}$ and hence on the degrees of the polynomials involved).

In order to do this we need again special polynomials which show us
the normalization of our isomorphism. We recall first how this
isomorphism is described:

Given $P_1^{(0)},P_2^{(0)} \in  U_\alpha^{(0)}$ we defined the
polynomial
$P_1^{(0)}\otimes P^{(0)}_2$ by
$(P_1\otimes P_2)(x)= \langle\langle P_1^{(0)}(d), P_2^{(0)}(xd \bar x) 
\rangle\rangle_0,$ where $\langle\langle \cdot,\cdot 
\rangle\rangle_0$ denotes the invariant scalar product chosen in
$U_\alpha^{(0)}.$ 

\medskip
We consider again the Gegenbauer polynomial $G^{(\alpha,0)}$ of
degree $\alpha$ in $U_\alpha^{(0)}\otimes U_\alpha^{(0)} $, derived in
the same way from the one-variable polynomial with indices $l=\alpha,
p=1/2$ given in \cite{vilenkin} as we did it above for the $G^{\beta}$
in $U_\beta$ and let  $\langle\langle \cdot,\cdot 
\rangle\rangle_0$ be normalized such that this polynomial is a
reproducing kernel, this normalization determines then our choice of
the isomorphism between $U_{2\alpha}$ and $U_\alpha^{(0)}\otimes
U_\alpha^{(0)}.$

\medskip
In order to relate the integrals in (\ref{eq:int_in4}) and
(\ref{eq:int_in3}) we evaluate them for a special choice of
polynomials:
We put $Q_1^{(0)}=Q_2^{(0)}=G^{(a'+b/2,0)}_z$ and
$Q_3^{(0)}=G^{(a',0)}_z$ for some quaternion $z$ of norm $1$ and trace $0$.

The integral in (\ref{eq:int_in3}) is then by \cite[p.490]{vilenkin}
equal to
\begin{eqnarray}
  \label{eq:c1}
  c_1&:=&\frac{\Gamma(\frac{3a'+b}{2}+1)\Gamma(\frac{a'+1}{2})^2
    \Gamma(\frac{a'+b+1}{2})}{\pi
    \Gamma(\frac{a'}{2}+1)^2\Gamma(\frac{a'+b}{2}+1)
    \Gamma(\frac{3a'+b+3}{2})}\\  
&=&16 \frac{\Gamma(a')^2\Gamma(a'+b)\Gamma(\frac{3a'+b}{2}+1)^2}
  {\Gamma(3a'+b+2)\Gamma(\frac{a'}{2})^2\Gamma(\frac{a'}{2}+1)^2
  \Gamma(\frac{a'+b}{2}) \Gamma(\frac{a'+b}{2}+1)}
\end{eqnarray}
(where the second form is derived from the first using the duplication
formula for the $\Gamma$-function and where a factor
$\Gamma(2w)/\Gamma(w)$ has to be replaced by $2$ if $w=0$).

\medskip
On the other hand, the definition of the isomorphism
 between $U_{2\alpha}$ and $U_\alpha^{(0)}\otimes
U_\alpha^{(0)}$ and the reproducing property of the Gegenbauer
polynomials imply that 
$$(G_z^{(\alpha,0)}\otimes G_z^{(\alpha,0)})(x)=G_z^{(\alpha,0)}(\bar
x z x),$$
and hence that 
\begin{multline}
  \label{eq:normal_1}
  \int (G_z^{((a'+b/2),0)}\otimes G_z^{((a'+b/2),0)})(x)
  (G_z^{((a'+b/2),0)}\otimes G_z^{((a'+b/2),0)})(x)\\
\times  (G_z^{(a',0)}\otimes G_z^{(a',0)})(x)dx\\
=\text{vol}(\text{Stab}(z))\int G_z^{((a'+b/2),0)}(z')G_z^{((a'+b/2),0)}(z')G_z^{(a',0)}(z')dz',
\end{multline}
where Stab$(z)$ is the set of $x$ of norm 1 with $\bar x z x=z.$

The normalizations of the integrals over the $3$-sphere and over the
$2$-sphere in \cite{vilenkin} are such that
$\text{vol}(\text{Stab}(z))=\frac{\pi}{4}$ holds.

Taken together we obtain 
\begin{multline}
  \label{eq:int_compare}
  \int
(Q_1^{(0)}\otimes Q_1^{(0)})(x)( Q_2^{(0)}\otimes
Q_2^{(0)})(x)(Q_3^{(0)}\otimes Q_3^{(0)})(x)dx\\
=\frac{\pi}{4c_1}
\left(\int
Q_1^{(0)}(z)Q_2^{(0)}(z)Q_3^{(0)}(z)dz\right)^2,
\end{multline}
where $c_1$ is the constant computed in (\ref{eq:c1}).

\medskip
This gives us the first formula for the central critical value of the
triple product $L$-function.

\bigskip
\begin{multline}
  \label{firstformula}
N^{-1} 2^{12a'+4b-4-\omega(N)}\pi^{5+6a'+2b}\la f_1,f_1\ran\la f_1,f_1\ran \la
f_3,f_3 \ran\\
\times{\frac {(b+a')a'
\Gamma
(\frac{a'+b}{2})^{2}\Gamma(\frac{a'}{2})^4\Gamma(3a'+b+2)}{
\Gamma(a')^{2}\Gamma(2a')\Gamma(\frac{3a'+b}{2}+1)^2\Gamma(a'+b)\Gamma(2a'+b+1)^2}}\\  
\times \bigl(
\int P_1(x)P_1(x)P_3(x)dx
\bigr)^2 
\end{multline}

In this we replace the Petersson product $\la f_i,f_i\ran$ by 
$$(4\pi)^{1-k_i}\Gamma(k_i)D_{f_i}(k-1)$$ (where $D_{f_i}$
denotes the symmetric square $L$-function of $f_i$), which
leads us to  
\begin{multline}
  \label{secondformula}
2^{-9-\omega(N)}\pi^{2}
D_{f_1}(k_1-1)D_{f_1}(k_1-1)D_{f_3}(k_3-1)\\   
\times{\frac {(a')^2(2a'+1)(a'+b)(b+2a'+1)^2
\Gamma
(\frac{a'+b}{2})^{2}\Gamma(3a'+b+2)\Gamma(\frac{a'}{2})^4} {\Gamma(a')^{2}\Gamma (\frac{3a'+b}{2}+1)^2\Gamma (a'+b)}}\\ 
\times \bigl(
\int P_1(x)P_1(x)P_3(x)dx
\bigr)^2 
\end{multline}

Here the factor 
$D_{f_1}(k_1-1)D_{f_1}(k_1-1)D_{f_3}(k_3-1)$ does not
contribute in an essential way
to the asymptotics as $k_1 \rightarrow \infty$
since it is well known that $k_i^{-\delta}<< D_{f_i}(k_i-1)<<k_i^\delta$
for all $\delta>0$, see e.g. \cite{Hoff-Lock}.

We analyze the total factor on the right hand side in front of 
$$D_{f_1}(k_1-1)D_{f_1}(k_1-1)D_{f_3}(k_3-1)\bigl(
\int P_1(x)P_2(x)P_3(x)dx
\bigr)^2 $$

with Stirling's formula:
For the first assertion of the proposition we fix $\nu_3$
and let $\nu_1$ tend to infinity; we  find
that the factor from above can for all $\epsilon>0$  be bounded from below by 
$$c'b^{3-\epsilon}$$ for some nonzero constant $c'$
depending on $a', \epsilon$ and the 
level $N$ as $b$ tends to infinity.

For the second part of the proposition we have all the $\nu_i$ equal,
which implies $b=0$ in the notation used above; we find that the factor
can (for all $\epsilon>0$) be bounded from below by $$c''
a^{5-\epsilon}$$
for some nonzero constant $c''$
depending on $\epsilon$ and the 
level $N.$
\medskip
We have to adjust a final normalization:
The $\vp_i$ were normalized to have $$\dla \vp_i,\vp_i\dra_0 =1$$
whereas we want them to have $L^2$-norm $1.$

A comparison of $\dla \vp_i,\vp_i\dra_0$ with the scalar product on 
the space of $\vp$ derived from the $L^2$-norm with the help of
\cite[p. 461]{vilenkin} shows that we have 
$$P_1=\sqrt{a'+(b+1)/2}^{-1}\,\,\tilde{P_1},\quad P_3=\sqrt{a'+1/2}^{-1}\,\,\tilde{P_3},$$
where the $\tilde{P_i}$ are $L^2$-normalized.

This multiplies the last formula with 
$$(a'+(b+1)/2)^{-2}(a'+1/2)^{-1}$$ in the first case and leads to an
expression that is bounded from below for every $\epsilon>0$ by 
$$C_1b^{1-\epsilon}$$ for some nonzero constant $C_1$ depending
on $\epsilon, a'$ and the 
level $N$ as $b$ tends to infinity.

In the second case the formula gets multiplied with
$(a'+1/2)^{-3}$  and leads to an
expression that is bounded from below for every $\epsilon>0$ by 
$$C_2a^{2-\epsilon}$$ for some nonzero constant $C_2$ depending
on $\epsilon$ and the 
level $N$ as $a$ tends to infinity.
This finishes the proof of the Proposition.

\bigskip
On the other hand one can investigate the dependence of the
central critical value $L(f_1,f_1,f_3;\frac{2k_1+k_3}{2}-1)$
on the level and weights with analytic methods from the
theory of $L$-functions.

We do this first for the equidistribution of mass conjecture, i.\ e.,
for the situation in which $\nu_3$ is fixed and $\nu_1=\nu_2$ tends to
infinity:

As usual in the theory of $L$-functions the first step is to
establish the convexity bound.

\begin{lemma}
Let $f_1,f_3$ be newforms of level $N$ and weights $k_1,k_3$
as above. Then

$$L(f_1,f_1,f_3;\frac{2k_1+k_3}{2}-1)=O_N(k_1^{1+\epsilon})$$
for all $\epsilon >0.$
\end{lemma}
{\it Proof.}
In \cite{IS} a general description of the convexity bound of (standard)
automorphic L-functions for $GL(n)$ is given. That bound is also applicable
for our triple $L$-function.\\
We recall the (normalized) functional equation of the triple $L$-function
(quoting from \cite{GK} for weight 2 and more generally from\cite{bs-triple} 
; we restrict ourselves
to the case where all the cusp forms involved are newforms of the same
-squarefree- level):
Putting
\begin{multline}
\Lambda(s)= \Gamma_{\mathbb C}(s+k_1+{k_3\over 2})
\Gamma_{\mathbb C}(s+1+{k_3\over 2})
\Gamma_{\mathbb C}(s+1+{k_3\over 2})\\
\times\Gamma_{\mathbb C}(s+1+k_1-{k_3\over 2})
{\mathcal L}(f_1,f_2,f_3,s)
\end{multline}
with
\begin{equation}\label{script-L}{\mathcal
  L}(f_1,f_2,f_3,s)=L(f_1,f_2,f_3,s+{k_1+k_2+k_3-3\over
  2})
\end{equation}
 and $\Gamma_{\mathbb C}(s)=(2\pi)^{-s}\Gamma(s)$
we get the functional equation
$$\Lambda(1-s)=(N^5)^{s-{1\over 2}}w \Lambda(s)$$ 
Under such circumstances, the convexity bound (as described in \cite{IS}) is
$${\mathcal L}(f_1,f_2,f_3,{1\over 2}+it)<<_{\epsilon}
(C(f_1,f_2,f_3,t))^{{1\over 4}+\epsilon}$$
Here $C(f_1,f_2,f_3,t)$ is given in terms of the gamma factors; 
in our case  this means
\begin{multline}C(f_1,f_2,f_3,t)=(N^5)
(1+\mid it -k_1-{k_3\over 2}\mid^2)\\ \times
(1+\mid it -1-{k_3\over 2}\mid^2)^2
(1+\mid it -1-k_1+{k_3\over 2}\mid^2)
\end{multline}
In our special case (i.e.$k_1=k_2$ and $k_3$ fixed) this implies an 
estimate of type
$${\mathcal L}(f_1,f_1,f_3,{1\over 2})<< k_1^{1+\epsilon}$$
or
$${\mathcal L}(f_1,f_1,f_3,{1\over 2})<< b^{1+\epsilon}$$

\medskip
For our intended application this result is just too weak.

We will therefore assume henceforth that one can break
convexity for the estimate of
$L(f_1,f_1,f_3;\frac{2k_1+k_3}{2}-1)$ in the $k_1$-aspect (see
\cite{IS} for a survey of subconvex estimates).

\medskip

{\bf Subconvexity hypothesis.} {\em Fix $f_3$ as above.
There is $\delta >0$ such that for all $f_1$ as above 
 $${\mathcal L}(f_1, f_1, f_3,\frac{1}{2}) =
 O_{f_3}(k_1^{1-\delta}).$$} 
\vskip0.5cm

Then the result from our proposition immediately translates into a 
statement about equidistribution of measures on the unit sphere that are
associated to Hecke eigenfunctions on a quaternion algebra which
proves the last assertion in the introduction (concerning Conjecture 1).

 \begin{proposition} \label{equi-propo}
Let $D$ be as above, let $S^2 \subseteq \Bbb R^3$ be identified with\\
$\{x \in D_{\infty}^{(0)}~|~n(x) = 1\}$ as above. For a harmonic
polynomial $P \in \Bbb C[X_1,X_2,X_3]$ let the measure $\mu_P$ on
$S^2$ be defined by
 $$\int_{S^2} f(x) d\mu_P(x) := \int_{S^2} |P(x)|^2 f(x) dx.$$
Then under our subconvexity hypothesis the measures $\mu_P$ become
equidistributed if $P$ runs through Hecke eigenfunctions of degree
$\nu$ for $\nu \longrightarrow \infty$, i.e., one has
 $$\lim_{\nu(P)\to\infty} \int_{S^2} f(x)d\mu_P(x) = \int_{S^2}f(x) dx
  \eqno{(\ast)}$$
for all continuous $f$ on $S^2$.
 \end{proposition}
{\it Proof.} We have to check $(\ast)$ only for Hecke eigenfunctions
$f = Q$, as these form a Hilbert space basis of $L^2(S^2)$.

Then for $Q \not= 1$ the right hand side of $(\ast)$ is zero,
for $Q = 1$ it is 1.

For $Q = 1$ we have equality in $(\ast)$ for all $P$ (of $L^2$-norm 1).

For $Q \not= 1$, Proposition 2 together with the convexity breaking
assumption implies that
 $$\lim_{\nu\to \infty} \int_{S^2} Q(x) |P(x)|^2 dx = 0,$$
which proves the assertion.

\medskip
We end with the statement and proof of Conjecture 1 for the case $q=3$

\begin{proposition}
Let $D$ and $P$ (harmonic of degree $\nu(P)$) be as in Proposition
\ref{equi-propo}. 
Then 
\begin{equation}
\lim_{\nu(P) \rightarrow \infty} \int_{S^2}P(x)^3dx=0
\end{equation}
\end{proposition}

{\it Proof.} According to (\ref{eq:mainresult1a}) we have for
$\epsilon>0$
\begin{equation}
  \vert \int_{S^2} P(x)^3 dx\vert^2 << \nu(P)^{\epsilon-2}L(f_1,f_2,f_3;2+3\nu(P))
\end{equation}
Now with the notation from (\ref{script-L}) (i.\ e.\ denoting by
$\mathcal L$ an $L$-function normalized to have functional equation
under $s \mapsto (1-s)$) we have 
\begin{equation}
{\mathcal L} (f,f,f,s)={\mathcal L}(\text{Sym}^3 f,s)({\mathcal L}(f,s))^2.
\end{equation} 
Furthermore $\mathcal L(\text{Sym}^3f,s)$ is the $L$-function of a
$GL_4$-cusp form \cite{kim-shahidi}, so we may apply the general
Molteni subconvexity bound (see \cite{IS}) to ${\mathcal
  L}(\text{Sym}^3f,\frac{1}{2}).$ This combined with the subconvex
bound for ${\mathcal L}(f,\frac{1}{2})$ due to Peng \cite{peng} shows
that $$\nu(P)^{\epsilon-2}L(f_1,f_2,f_3;2+3\nu(P))=O(\nu^{-\delta})$$
for a fixed $\delta>0.$ This proves the Proposition.

\medskip
{\it Remark.} We see presently no way to extend the statement of
Proposition \ref{firstproposition} (and hence our arguments in this
article) to a product of more than three
polynomials, since our proofs here and in \cite{bs-triple} use several
special features of the case of three polynomials. In particular 
our proofs depend on 
\begin{itemize}
\item the  existence of an
integral representation for the triple product $L$-function using
an Eisenstein series whose special value is expressed by theta series
with spherical harmonics 
\item the existence and uniqueness of
trilinear forms 
on tensor products of spaces of harmonic polynomials (and their
explicit description by Ibukiyama's 
generating series in \cite{ibu}).
\end{itemize}

\medskip
{\em Remark.} In the case that the class number  of the quaternion
algebra is  $h\ne 1$ Hecke eigenforms (on the adelic
quaternion algebra) give rise to $h$-tuples of harmonic
polynomials, they should then be viewed
as functions on the disjoint sum of $h$ copies of the unit sphere. All
arguments from above can be carried out for such $h$-tuples
of harmonic polynomials (resp.\ functions on the disjoint sum of $h$
copies of the unit sphere).

\vskip0.5cm
{Siegfried B\"ocherer,
Kunzenhof 4B,
79117 Freiburg,
Germany,\\
boech@siegel.math.uni-mannheim.de\\

Peter Sarnak,
Department of Mathematics,
Princeton University,
Fine Hall, Washington Road,
Princeton, NJ 08544,
sarnak@math.princeton.edu}\\

Rainer Schulze-Pillot,
Fachrichtung 6.1 Mathematik,
Universit\"at des Saarlandes (Geb. 27.1),
Postfach 151150,
66041 Saarbr\"ucken,
Germany,
schulzep@math.uni-sb.de

\newpage
\begin{center} 
{\bf {\large Appendix} \\ 
\vskip0.5cm S. B\"ocherer, R. Schulze-Pillot: Corrections to our
  article  ``On
  the central 
critical value of the triple product L-function'' (reference
  \cite{bs-triple} of this article)}
\end{center}\vskip1cm

\begin{itemize}
\item {\bf p.5:} \\
l.8 f.a.: read $\left(L_{\alpha}^{(b)}f\right)$ rather than 
$\left(L_{\alpha}^{(b)}\right)$
\item {\bf p.5:}\\
l.9 f.b. read $\partial_{12}X_2$ rather than $\partial_{12}X_1$ 
\item {\bf p.7:}\\
l.12 f.b.: read "(1.9)" rather than "(1.10)"
\item {\bf p.7:} \\
Skip the sentences "It is easy..." (l.8 f.b.) until "...polynomial in
$\partial_{12},\partial_{13},\partial_{23}$ (l.5 f.b.)
\item {\bf p.7:}\\
l.4.f.b.: read $a'$ rather than $a$
\item {\bf p.7:} \\
l.3 f.b.:Formula (1.11) should read
$${2^b 2^{a'}(2i)^{3a'}\over
(\alpha+a')^{[b]} b!}
\left(\partial_{12}\partial_{13}
\partial_{23}\right)^{a'}\left(\partial_{12}X_2+\partial_{13}X_3\right)^b$$
\item {\bf p.9:}\\
l.11 f.a.: read $a'$ rather than $a$  
\item {\bf p.12}\\
l.8 f.b.: the formula should be
$$A(s,b)=(-1)^b\cdot
{ 2^{-a}(2s-2)^{[a]}(s)^{[a']}(2s+a-2)^{[b]}\over
b!(s+a'-1)^{[b]}
}$$
\item {\bf p.13:} \\
in formula (2.11)
$$\left(g<Z>^*\right)  \qquad\mbox{rather than} \quad \times <g<Z>^*>$$
$${\det}^{k+a}\qquad \mbox{rather than}\quad \det^{k+a}$$
\item {\bf p.17:}\\
l.7 f.b.: read $h'$ rather than $h$
\item {\bf p.18:}\\
delete factor $2$ in formula (2.28)
\item {\bf p.32}\\
in formula (4.1) read $L_p(f,\phi,\psi,s+3a'+2r+b-2)$
\item {\bf p.33}\\
Formula (4.3) is correct only under the additional assumption $r=2$ 
\item {\bf p.33}\\
In (4.3) read $\zeta$ rather than zeta
\item {\bf p.33}\\
In the functional equation (4.4) the exponent of $N$ should be\\
$-4(s-{k_1+k_2+k_3\over 2}+1)$;
the exponent of $gcd(N_f,N_{\phi},N_{\psi})$
should be\\
$-(s-{k_1+k_2+k_3\over 2}+1)$
\item {\bf p.34:}\\
l.10: read allows for $p \mid N$ rather than allows for $p \nmid N$ 
\item {\bf p.35:}\\
l.5.f.a.: read (2.1) rather than (2.2)
\item {\bf p.37}:\\
l.9 f.b. read $i^{3a'}\pi^{3a'+b}$ rather than $\pi^{3a+2b}$
\item {\bf p.37:}\\
l.7 f.b: read $\left(i^{3a'}\pi^{3a'+b}\right)^{-1}$ rather than 
$\pi^{-3a-2b}$ 
\item {\bf p.38:}\\
in formula 5.4: read $i^{3a'}\pi^{3a'+b}$ rather than 
$\pi^{3a+2b}$ 
\item {\bf p.42:}\\
In Lemma 5.5 and in line 19 read
$$\left(\Sum_{i=1}^h\frac{T_0(\vp_f(y_i)\otimes
\vp_\phi(y_i) \otimes\vp_\psi(y_i))}{e_i}\right)$$
\item {\bf p.44:}\\
The first line of formula (5.9) should read
\begin{eqnarray*}
(-1)^{\omega(N)+\omega(M_1,M_2)}
2^{5+3b+8a'-\omega(gcd(N_f,N_{\phi},N_{\psi})}\pi^{5+6a'+2b}\\
\times N^{2}(M_1M_2)^{-3}M_3^{-6}
{1\over {b\choose \nu_2}}{ b!(a'+1)^{[b]}\over
2^{[a]}2^{[a']}(a+2)^{[b]}}\quad\quad\\
\times 
{(a'+b+1)\Gamma(2a'+b+2)\over
\Gamma(3a'+b+2)\Gamma(a'+\nu_2+1)\Gamma(a'+\nu_3+1)}
\end{eqnarray*}
\item {\bf p.45:}\\
l.8.f.a.: Nagoya Math.J.147(1997), 71-106
\item {\bf p.45} \\
l.3 f.b.: Comm.Math.Univ.S.Pauli 48(1999), 103-118
\end{itemize}
\end{document}